\documentclass{amsart}
\usepackage{amsmath}
\usepackage[ngerman,french,english]{babel}

\newtheorem{ques}{\textbf{Question}}

\newtheorem{theorem}{\textbf{Theorem}}

\newtheorem{corollary}{\textbf{Corollary}}

\def\Z {\mathbb{Z}}

\def\Q {\mathbb{Q}}

\def\C {\mathbb{C}}

\def\den{\mathop{\rm den}}

\theoremstyle{remark}

\numberwithin{equation}{section}

\begin{document}

\title[TRANSCENDENTAL ANALYTIC FUNCTIONS MAPPING $\mathbb{Q}$ INTO ITSELF]{A NOTE ON TRANSCENDENTAL ANALYTIC FUNCTIONS WITH RATIONAL COEFFICIENTS MAPPING $\mathbb{Q}$ INTO ITSELF}

\author[J. LELIS]{JEAN LELIS}
\address{FACULDADE DE MATEMÁTICA, ICEN, UNIVERSIDADE FEDERAL DO PARÁ, BELÉM, PA, BRAZIL}
\email{jeanlelis@ufpa.br}

\author[D. MARQUES]{DIEGO MARQUES*}
\address{DEPARTAMENTO DE MATEM\'{A}TICA, UNIVERSIDADE DE BRAS\'ILIA, BRAS\'ILIA, DF, BRAZIL}
\email{diego@mat.unb.br}

\author[C. G. MOREIRA]{CARLOS GUSTAVO MOREIRA}
\address{INSTITUTO DE MATEM\' ATICA PURA E APLICADA, RIO DE JANEIRO, RJ, BRAZIL}
\email{gugu@impa.br}

\author[P. TROJOVSK\' Y]{PAVEL TROJOVSK\' Y}
\address{FACULTY OF SCIENCE, UNIVERSITY OF HRADEC KR\' ALOV\' E,
CZECH REPUBLIC}
\email{pavel.trojovsky@uhk.cz}


%    General info


\begin{abstract}
In this note, the main focus is on a question about transcendental entire functions mapping $\mathbb{Q}$ into $\mathbb{Q}$ (which is related to a Mahler's problem). In particular, we prove that, for any $t>0$, there is no a transcendental entire function $f\in\mathbb{Q}[[z]]$ such that $f(\mathbb{Q})\subseteq\mathbb{Q}$ and whose denominator of $f(p/q)$ is $O(q^{t})$, for all rational numbers $p/q$, with $q$ sufficiently large.
\end{abstract}

\maketitle

%%  SECTION 1
\section{Introduction}

Transcendental number theory began in 1844, when Liouville \cite{liouville} proved the existence of transcendental numbers. In fact, he was able to explicit an infinite class of such numbers. These numbers are the well-known {\it Liouville numbers}: a real number $\xi$ is called a Liouville number, if there exists a sequence of distinct rational numbers $(p_k/q_k)_k$, with $q_k>1$, such that
\[
0<\left|\xi-\frac{p_k}{q_k}\right|<\frac{1}{q_k^{\omega_k}},
\] 
where $\omega_k$ tends to infinite as $k\to \infty$. The set of all Liouville numbers is denoted by $\mathbb{L}$.

In his pioneering book, Maillet \cite{mai}, in 1906, proved that $f(\mathbb
L)\subseteq \mathbb{L}$, for any non-constant rational function $f\in \mathbb{Q}(z)$. In light of this fact, in 1984, Mahler \cite{mahler1} raised the following question:
\begin{ques}\label{qm}
Are there transcendental entire functions $f(z)$ such that if $\xi$ is any Liouville number, then so is $f(\xi)$?
\end{ques}

In 2015, Marques and Moreira \cite{Mar} showed the existence of uncountably many transcendental entire functions $f$ such that $f(\Q)\subseteq \Q$ and for which $\den(p/q)<q^{8q^2}$, for all rational number $p/q$, with $q>1$ (here, and in what follows, $\den(z)$ denotes the denominator of the irreducible rational number $z$). It follows from their argument that Question \ref{qm} has a positive answer if the following question also has:

\begin{ques}\label{Q1}
Are there transcendental entire functions $f(z)$ such that $f(\mathbb{Q})\subseteq\mathbb{Q}$ and 
\[
 \den(f\left(p/q\right))= O(q^{t}),
 \]
for all rational numbers $p/q$, with $q$ sufficiently large, where $t\geq 0$ is a given positive integer?
\end{ques} 

There are some progress in this question. For instance, in 2016, Marques, Ramirez and Silva \cite{Mar1} showed that there is `no'\ transcendental entire function $f(z)\in\mathbb{Q}[[z]]$ such that $f(\mathbb{Q})\subseteq\mathbb{Q}$ and $\den f(p/q)=o(q)$ (this is the case $t\in [0,1)$ in Question \ref{Q1}). By using Whittaker's theory \cite{Whit1} of polynomial expansions of analytic functions, in 2020, Lelis and Marques \cite{LM} proved that the answer is also `no'\ for the case in which $f(z)\in\C[[z]]$ and $t=1$. 

The goal of this note is to prove, in particular, the non-existence of functions as in Question \ref{Q1} with rational coefficients. More precisely,

\begin{theorem}\label{main1}
	Let $t$ be a positive integer and $\Omega\subseteq \C$ be a neighborhood of origin. Then, there is no a transcendental function $f(z)=\sum_{k\geq 0}a_kz^k\in\mathbb{C}[[z]]$, analytic in $\Omega$, such that $a_k\in \Q$, for all $k\in [0, t]$, $f(\mathbb{Q}\cap\Omega)\subseteq\mathbb{Q}$ and
	\[
	\den f(p/q)=O(q^{t/2}),
	\]
	for all rational numbers $p/q\in\Omega$, with $q$ sufficiently large.
\end{theorem}

As an immediate consequence, we infer that
\begin{corollary}\label{cor1}
	Let $t$ be a positive integer. Then, there is no a transcendental entire function $f(z)\in\mathbb{Q}[[z]]$ such that  $f(\mathbb{Q})\subseteq\mathbb{Q}$ and
	\[
	\den f(p/q)=O(q^t),
	\]
	for all rational numbers $p/q$, with $q$ sufficiently large.
\end{corollary}

\section{The Proof of Theorem \ref{main1}}

In order to simplify the argument, and causing no loss of generality, we shall prove the theorem for $2t$ (instead of $t$). Aiming for a contradiction, let $t\geq 1$ be an integer and suppose that $f(z)\in\C[[z]]$ is a transcendental function which is analytic in the neighborhood $\Omega\subseteq \C$ of the origin given by
\[
f(z)=\sum_{k=0}^{\infty}a_kz^k
\]
such that $a_k\in\Q$, for all $ k\in[0,2t]$, $f(\mathbb{Q}\cap\Omega)\subseteq\mathbb{Q}$ and $\den(f(p/q))=O(q^t)$, for all $p/q\in\Omega\cap \Q$ with $q$ sufficiently large. Note that, without loss of generality, we may assume that $f(0)=0$, that is, $a_0=0$. 

Let $F(z)\in\Q[z]$ be the polynomial given by
\[
F(z)=\sum_{n=1}^{2t}a_nz^n.
\]
Then, there exist polynomials $P(z)$ and $Q(z)$ in $\Q[z]$ of degrees at most $t$ such that
\[
R(z):=\frac{P(z)}{Q(z)}=F(z)+z^{2t+1}K(z),
\]
where $K(z)\in\C[[z]]$ is analytic in a neighborhood of origin. Indeed, suppose that
\[
Q(z)=q_0+q_1z+\cdots+q_tz^t,
\]
and consider the product $Q(z)F(z)$ given by
\[
Q(z)F(z)=b_1z+\cdots+b_{3t}z^{3t}.
\]
Therefore, we need to determine $q_0,q_1,\ldots,q_t$ in $\Q$, not all zero, such that $b_j=0$, for all $j\in [t+1, 2t]$. In other words, we want a non-trivial rational solution for the $t\times(t+1)$ homogeneous linear system
\begin{eqnarray*}
	\left\{\begin{array}{ccccccccc}
		a_{2t}q_0 & + & a_{2t-1}q_1 & + & \cdots & + & a_{t}q_t & = & 0 \\
		a_{2t-1}q_0 & + & a_{2t-2}q_1 & + & \cdots & + & a_{t-1}q_t & = & 0 \\
		\vdots & \vdots & \vdots & \vdots & \vdots & \ddots & \vdots & \vdots \\
		a_{t+1}q_0 & + & a_{t}q_1 & + & \cdots & + & a_1q_t & = & 0
	\end{array}\right.
\end{eqnarray*}
in the variables $q_0,q_1,\ldots,q_t$. Since the number of variables is larger than the number of equations, then a basic Linear Algebra result ensures the existence of a solution $(q_0,\ldots,q_t)\in \Q^{t+1}\backslash \{(0,\ldots, 0)\}$. Now, let $j\geq 0$ be the smallest non-negative integer such that $q_j\neq 0$. Thus, we may suppose that $q_j=1$ (because $(0,\ldots, 0, 1,q_{j+1}/q_j,\ldots, q_t/q_j)$ is also a solution of the previous linear system). So, we define
\[
P(z):=b_{j+1}z^{j+1}+b_{j+2}z^{j+2}+\cdots+b_tz^t.
\]
Therefore
\[
P(z)=Q(z)F(z)+z^{2t+1}S(z),
\]
with $S(z)\in\Q[z]$ and so
\begin{equation}\label{eqRt}
   R(t)=\frac{P(z)}{Q(z)}=F(z)+z^{2t+1-j}K(z),
\end{equation}
where $K(z)=S(z)/(1+q_{j+1}z+\cdots+q_tz^{t-j})$ is analytic in a neighborhood $\Omega'\subseteq\C$ of origin. 

On the other hand, we have that 
\begin{equation}\label{eqf}
f(z)=F(z)+z^{2t+1}\Tilde{K}(z),
\end{equation}
where $\Tilde{K}(z)$ is analytic in $\Omega$. Let $M'$ be the smallest positive integer such that $1/M'\in \Omega\cap \Omega'$. Thus, combining \eqref{eqRt} together with \eqref{eqf} and for all positive integers $M>M'$, we get
\begin{equation}\label{mod1}
\left|f\left(\frac{1}{M}\right)-R\left(\frac{1}{M}\right)\right|=\frac{1}{M^{2t+1-j}}\Theta(M),
\end{equation}
where
\[
\Theta(M)=\left|\frac{1}{M^j}\Tilde{K}\left(\frac{1}{M}\right)-K\left(\frac{1}{M}\right)\right|.
\]

Now, by hypothesis, $f(1/M)=A_M/B_M\in\Q$, with $\gcd(A_M,B_M)=1$ and $B_M\leq C_1M^t$ (for some positive constant $C_1$). Moreover, we have that
\[
R\left(\frac{1}{M}\right)=\frac{b_{j+1}M^{t-j-1}+b_{j+2}M^{t-j-2}\cdots+b_t}{M^{t-j}+q_{j+1}M^{t-j-1}+\cdots+q_t}
\] 
is a rational number. Note that, if we suppose that $f(1/M)\neq R_t(1/M)$ for infinitely many integers $M>M'$, then 
\begin{equation}\label{mod2}
\left|f\left(\frac{1}{M}\right)-R\left(\frac{1}{M}\right)\right|\geq \frac{1}{C_1M^t(M^{t-j}+q_{j+1}M^{t-j-1}+\cdots+q_t)}\geq \frac{1}{C_2M^{2t-j}},
\end{equation}
for all integer $M$ sufficiently large (say $M>M_0>M'$), where $C_2=2C_1$. By \eqref{mod1} and \eqref{mod2}, we obtain that 
\[
\frac{1}{M^{2t-j}}\ll \left|f\left(\frac{1}{M}\right)-R\left(\frac{1}{M}\right)\right|=\frac{1}{M^{2t+1-j}}\Theta(M)\ll \frac{1}{M^{2t+1-j}}
\]
yielding to the absurdity that $M^{2t+1-j}\ll M^{2t-j}$, for infinitely many positive integers $M$ for which $f(1/M)\neq R(1/M)$.

Therefore, $f(1/M)=R(1/M)$ for all large enough integer $M$, say $M\geq M_1$. Thus, $f(z)$ and $R(z)$ coincide in the set $\{1/M: M>\max\{M_0,M_1\}\}\subseteq \Omega\cap \Omega'$ which has a limit point. By the Identity Theorem for Analytic Functions, we infer that $f(z)=R(z)$, for all $z\in\Omega\cap\Omega'$ with contradicts the transcendence of $f(z)$. The proof is then complete.
\qed
%Acknowledgements

\section*{Acknowledgement}
D.M was supported by CNPq-Brazil. Part of this work was done during a visit of D.M and J.L to IMPA (Rio de Janeiro) which provided them excellent working conditions. 
P.T. was supported by the Project of Excellence, Faculty of Science, University of Hradec Kr\'alov\'e, No. 2210/2023-2024.

% \appendix

% \section{}
% \label{}

%% BIBLIOGRAPHY

\end{document}